\documentclass[11pt]{amsart}
\newtheorem{theorem}{Theorem}[section]

\newtheorem{lemma}[theorem]{Lemma}
\newtheorem{proposition}[theorem]{Proposition}

\theoremstyle{definition}
\newtheorem{definition}[theorem]{Definition}
\theoremstyle{remark}

\numberwithin{equation}{section}

\begin{document}

\title[A simple proof of the optimal power in Liouville theorems]
{A simple proof of the optimal power in Liouville theorems}
\author{Salvador Villegas}
\thanks{The author has been supported by Ministerio de Ciencia, Innovaci\'on y Universidades of Spain PGC2018-096422-B-I00 and by Junta de Andaluc\'{\i}a  A-FQM-187-UGR18.}
\address{Departamento de An\'{a}lisis
Matem\'{a}tico, Universidad de Granada, 18071 Granada, Spain.}
\email{svillega@ugr.es}

\begin{abstract}
Consider the equation div$(\varphi^2 \nabla \sigma)=0$ in $\mathbb{R}^N,$ where $\varphi>0$. It is well-known \cite{BCN, AC} that if there  exists $C>0$ such that  $\int_{B_R}(\varphi \sigma)^2 dx\leq CR^2$ for every $R\geq 1$, then $\sigma$  is necessarily constant. In this paper we present a simple proof that this result is not true if we replace $R^2$ with $R^k$ for $k>2$ in any dimension $N$. This question is related to a conjecture by De Giorgi \cite{DG}. 

\end{abstract}

\maketitle

{\footnotesize \noindent 2010 Mathematics Subject Classification: 35B08, 35B35, 35J91. 

}

\section{Introduction and main results}

In 1978, E. De Giorgi \cite{DG} stated the following conjecture:

{\bf Conjecture.} Let $u: \mathbb{R}^N \rightarrow (-1,1)$ be a smooth entire solution of the Allen-Cahn equation $-\Delta u=u-u^3$ which is monotone in one direction (for instance $\partial u/\partial_{x_N}>0$ in $\mathbb{R}^N$). Then $u$ depends only on one variable (equivalently, all its level sets are hyperplanes), at least if $N\leq 8$.

This conjecture was proved in 1997 for $N=2$ by Ghoussoub and Gui \cite{GG}, and in 2000 for $N=3$ by Ambrosio and Cabr\'e \cite{AC} . In 2011, del Pino, Kowalczyk, and Wei \cite{DKW} established that the conjecture does not hold for $N\geq 9$, as suggested in De Giorgi's original statement. In dimensions $4\leq N\leq 8$ the conjecture remains still open.

It is easily seen that a monotone solution of the Allen-Cahn equation is stable in the following sense:

\begin{definition}\label{stability}
Let $G\in C^2(\mathbb{R})$. We say that a solution $u\in C^2(\mathbb{R}^N)$ of $\Delta u=G'(u)$ in $\mathbb{R}^N$ is stable if
$$Q(v):=\int_{\mathbb{R}^N}\left( \vert \nabla v \vert^2 +G''(u)v^2\right)  dx\geq 0,$$
\noindent for every $v\in C^1(\mathbb{R}^N)$ with compact support in $\mathbb{R}^N$.

\end{definition}

Note that the above expression is nothing but the second variation of the energy functional in a bounded domain $\Omega\subset \mathbb{R}^N$:

$$E_\Omega (u)=\int_\Omega\left( \frac{1}{2}\vert \nabla u \vert^2 +G(u)\right) dx.$$

On the other hand, note that the Allen-Cahn equation $-\Delta u=u-u^3$ corresponds to the function $G(s)=(1-s^2)^2/4$.

To prove the conjecture for $N\leq 3$, it is used the notion of the stability and the following Liouville-type theorem due to Ambrosio and Cabr\'e \cite{AC}, which was motivated by a simpler version in \cite{BCN}.

\begin{theorem}{\bf (Ambrosio-Cabr\'e \cite{AC}).}\label{Lio}
Let $\varphi\in L_{loc}^\infty (\mathbb{R}^N)$ be a positive function. Assume that $\sigma\in H_{loc}^1 (\mathbb{R}^N)$ satisfies

\begin{equation}\label{ine}
\sigma \, \mbox{div}\, (\varphi^2\nabla\sigma )\geq 0 \ \ \ \mbox{in }\mathbb{R}^N
\end{equation}

\noindent in the distributional sense. For every $R>0$, let $B_R=\{ \vert x\vert<R\}$ and assume that there exists a constant independent of $R$ such that

$$
\int_{B_R}(\varphi \sigma)^2 dx\leq CR^2  \ \ \ \mbox{for every }\ R\geq 1.
$$

Then $\sigma$ is constant.

\end{theorem}

The idea is the following. If $u$ is a solution in De Giorgi's conjecture, consider the functions $\varphi:=\partial u/\partial_{x_N}>0$ and $\sigma_i:=\partial_{x_i}u/ \partial_{x_N}u$, for $i=1,...,N-1$.  Since both $\partial_{x_i}u$ and $\varphi$ solves the same linear equation $-\Delta v=(1-3u^2)v$, an easy computation shows that div$(\varphi^2\nabla\sigma_i )=0$. In dimensions $N\leq 3$ it is proved that $\int_{B_R}\vert\nabla u\vert^2 dx\leq CR^2$, for every $R\geq 1$. Applying Theorem \ref{Lio} gives $\sigma_i$ is constant and it follows easily that $u$ is a one-dimensional function. Observe that in this reasoning it is just used div$(\varphi^2\nabla\sigma_i )=0$, which is a an stronger condition than (\ref{ine}).

Motivated  by the useful applications of  Liouville-type theorems to these kind of problems, many authors have posed questions that allow us to understand qualitative properties of the solutions of general nonlinear partial differential equations of the form $\Delta u=G'(u)$ in $\mathbb{R}^N$. Here we show some of them.

{\bf Question A (Berestycki-Caffarelli-Nirenberg \cite{BCN}).} Let $L=-\Delta-V$ be a Schr{\"o}dinger operator on $\mathbb{R}^N$ with a smooth and bounded potential $V$. Suppose that $u\in W^{2,p}$ for some $p>N$ is a bounded and sign-changing solution for $Lu=0$. Set 
$$\lambda_1(V)=\inf\left\{\frac{\int_{\mathbb{R}^N}\vert\nabla\psi\vert^2-V\psi^2}{\int_{\mathbb{R}^N}\psi^2}; \psi\in V^\infty_0({\mathbb{R}^N})\right\}.$$
Then is $\lambda_1(V)<0$?

{\bf Question B (Berestycki-Caffarelli-Nirenberg \cite{BCN}).} Let $\varphi\in L_{loc}^\infty (\mathbb{R}^N)$ be a positive function. Assume that $\sigma\in H_{loc}^1 (\mathbb{R}^N)$ is a weak solution of   div$(\varphi^2 \nabla \sigma)=0$. If $\varphi\sigma$ is bounded, then is $\sigma$ constant?

In \cite{BCN} a positive answer to Question A is deduced from a positive answer to Question B. By Theorem \ref{Lio}, it is easily deduced that the answer to Question B is ''yes'' in dimensions $N=1,2$ and therefore also to Question A.  In \cite{GG}  Ghoussoub and Gui proved that the answer to Question A is ''no'' if $N\geq 7$  (and therefore also to Question B) by exploiting the idea of differentiating with respect to the first variable $x_1$ any solution of the semilinear pde $\Delta v+v^p=0$, $v>0$, $v(x)\to 0$ when $\vert x\vert \to\infty$ in the whole space $\mathbb{R}^N$, for an appropriate exponent $p>1$. A short time later, Barlow \cite{B} used probabilistic methods to construct counterexamples giving a negative answer to Questions A and B in any dimension $N\geq 3$.

Another question related to the previous ones, which in fact is the main point of this paper, is raised by Alberti, Ambrosio and Cabr\'e:

{\bf Question C (Alberti-Ambrosio-Cabr\'e \cite{AAC}).} If $0<\varphi\in L_{loc}^\infty (\mathbb{R}^N)$ and $\sigma\in H_{loc}^1 (\mathbb{R}^N)$ satisfy

\begin{equation}\label{equ}
\mbox{div}\, (\varphi^2\nabla\sigma )=0 \ \ \ \mbox{in }\mathbb{R}^N
\end{equation}

\noindent in the distributional sense, what is the optimal (maximal) exponent $\gamma_N$ such that

$$\int_{B_R}(\varphi \sigma)^2 dx\leq CR^{\gamma_N}  \ \ \ \mbox{for every }\ R\geq 1 \ \Longrightarrow \sigma \mbox{ constant}?$$

By the result of Barlow \cite{B}, it is deduced that $\gamma_N<N$ when $N\geq 3$. Also, a sharp choise in the counterexamples of Ghoussoub and Gui \cite{GG} shows that $\gamma_N<2+2\sqrt{N-1}$ for $N\geq 7$. More recently, Moradifam \cite{M} used some ideas of the probabilistic methods of Barlow \cite{B} to prove that  $\gamma_N<3$ when $N\geq 4$. 

On the other hand, assuming only inequality (\ref{ine}) (instead of equality (\ref{equ}))  Gazzola \cite{GA} proved the sharpness of the exponent $2$ in Theorem \ref{Lio} using nondifferentiable counterexamples. In any case, we would like to emphasize that we are more interested in the case of equality (\ref{equ}) than of inequality (\ref{ine}), since in the study of stable solutions it is precisely equality that is obtained.

In this paper we present a simple proof that $\gamma_N=2$ for every $N\geq1$. In other words, the exponent $2$ in Theorem \ref{Lio} is sharp, when equality  holds in (\ref{ine}). 

In fact, our obtained result could be deduced from a recent paper by the author \cite{yo}, in which more general functions (not only of the type $\Psi(R)=CR^k$) are considered in this kind of Liouville theorems. However, due to the simplicity of the proof of Theorem \ref{optimal} below and the original motivation of this type of questions, we think that it is appropriate to write the proof separately for the case of pure powers.

\begin{theorem}\label{optimal}

Let $k>2$ be a real number and $N\geq1$ any dimension. Then there exist $C>0$, $0<\varphi \in C^\infty (\mathbb{R}^N)$ such that equation  div $(\varphi^2 \nabla \sigma)=0$  admits a nonconstant classical solution $\sigma \in C^\infty (\mathbb{R}^N)$ satisfying

$$  \int_{B_R}(\varphi \sigma)^2 dx\leq CR^k \mbox{ for every }R\geq 1.$$ 

\end{theorem}

Therefore, to apply the same arguments as in the proof of the conjecture of De Giorgi for dimensions $N\leq 3$, it seems essential to establish a connection between $R^{\gamma_N}$ and $\int_{B_R}\vert\nabla u\vert^2$ (for varying $R$). Since $\vert \nabla u\vert$ is bounded we have the upper bound  $\int_{B_R}\vert\nabla u\vert^2\leq KR^N$ ($R>0$), for some $K>0$. In the following proposition we obtain a lower bound of this expression. Note that this lower bound cannot be improved, since in the case in which $u$ is one-dimensional we have $\int_{B_R}\vert\nabla u\vert^2\sim CR^{N-1}$ for large $R$. 

\begin{proposition}\label{porlomenos}
Let $G\in C^2(\mathbb{R})$ be a nonnegative function and $u\in C^2(\mathbb{R}^N)$ be a stable nonconstant bounded entire solution of $\Delta u=G'(u)$. Suppose that $G$ satisfies

\noindent  {\bf (H)}\ \ \ \ There exists $K>0$ such that $-(\sqrt{G})''(s)\geq K $ for every $s\in u(\mathbb{R}^N).$

 Then there exist $C,R_0>0$ such that

\begin{equation}\label{lowerbound}
\int_{B_R}\vert \nabla u\vert^2\geq CR^{N-1},\ \ \ \mbox{ for all }R\geq R_0.
\end{equation}

\end{proposition}

Note that in the classical Allen-Cahn equation $-\Delta u=u-u^3$, $\vert u\vert <1$, we have that $G(s)=(1-s^2)^2/4$ satisfies the hypothesis (H), since $-(\sqrt{G})''(s)=1$ for every $s\in(-1,1)$.

The lower bound (\ref{lowerbound}) follows also from a result of the forthcoming paper by Cabr\'e, Cinti and Serra \cite{CCS}, carried out independently of ours, which shows that Dirichlet energy in a ball controls potential energy in a slightly smaller ball. This bound together with Modica monotonicity formula (see Theorem \ref{modica2}) leads immediately to (\ref{lowerbound}).

Finally, to obtain a more precise study on the functional energy $E_{B_R}$, we provide the following result, which establishes that the Dirichlet and potential energies have the same behavior in $B_R$, for large $R$.

\begin{proposition}\label{energiasiguales}
Let $G\in C^2(\mathbb{R})$ be a nonnegative function and $u\in C^2(\mathbb{R}^N)$ be a stable nonconstant bounded entire solution of $\Delta u=G'(u)$. Suppose that $G$ satisfies (H). Then
$$\lim_{R\to\infty}\frac{\displaystyle{\int_{B_R}\frac{1}{2}\vert\nabla u\vert^2}}{\displaystyle{\int_{B_R}G(u)}}=1.$$
\end{proposition}

\section{Relation between Dirichlet and potential energies} 

In this section we prove some results concerning  the Dirichlet and potential energies. We will use the following two results due to Modica.

\begin{theorem}{\bf (Modica \cite{Mod1}).}\label{modica}
Let $G\in C^2(\mathbb{R})$ be a nonnegative
function and $u$ be a bounded solution of $\Delta
u=G'(u)$ in $\mathbb{R}^N$. Then,
$$
\frac{\vert\nabla u\vert^2}{2}\leq G(u) \quad {\rm in}\ \mathbb{R}^N.
$$
In addition, if $u$ is not constant, then 
$G(u(x))>0$ for all  $x\in\mathbb{R}^N$.
\end{theorem}

In \cite{Mod1} this bound was proved under the hypothesis
$u\in C^3(\mathbb{R}^N)$. The result as stated above, which applies to all
solutions ---recall that every solution is $C^{2,\alpha}(\mathbb{R}^N)$
since $G\in C^2(\mathbb{R})$--- was established by Caffarelli, Garofalo and Seg{\`a}la \cite{cafgarseg}.

\begin{theorem}{\bf (Modica \cite{Mod2}).}\label{modica2}
Let $G\in C^2(\mathbb{R})$ be a nonnegative
function and $u$ be a bounded solution of $\Delta
u=G'(u)$ in $\mathbb{R}^N$. Then,
$$
\Phi(R):=R^{1-N}\int_{B_R}\left( \frac{1}{2}\vert \nabla u \vert^2 +G(u)\right) dx
$$

\noindent is nondecreasing in $R\in (0,\infty)$. In particular, if $u$ is not constant then there exists a positive constant $c_0$ such that
$$\int_{B_R}\left( \frac{1}{2}\vert \nabla u \vert^2 +G(u)\right) dx\geq c_0 R^{N-1},\ \ \ \mbox{ for all }R\geq 1.$$
\end{theorem}

\begin{lemma}\label{caract}

Let $G\in C^2(\mathbb{R})$ be a nonnegative function and $u\in C^2(\mathbb{R}^N)$ be a nonconstant bounded entire solution of $\Delta u=G'(u)$. Then $u$ is stable if and only if
$$\int_{\mathbb{R}^N}G(u)\vert \nabla \mu \vert^2\, dx\geq \int_{\mathbb{R}^N}\left( G(u)-\frac{1}{2}\vert \nabla u\vert^2 \right) \left( -2(\sqrt{G})''(u)\right) \sqrt{G(u)}\mu^2 \, dx,$$
\noindent  for every $\mu\in C^1(\mathbb{R}^N)$ with compact support in $\mathbb{R}^N$.

\end{lemma}

\noindent\textbf{Proof.}

First of all note that, by Theorem \ref{modica}, we have $G(u(x))>0$ for every $x\in \mathbb{R}^N$. This means that $G(s)>0$ for every $s\in I:=u(\mathbb{R}^N)$. Then $\sqrt{G}\in C^2(I)$ and the above inequality has sense.

Let $\mu\in C^1(\mathbb{R}^N)$ have compact support and $\omega\in C^2(\mathbb{R}^N)$ . We obtain that

$$Q(\omega\mu)=\int_{\mathbb{R}^N}\left( \omega^2\vert\nabla\mu\vert^2+\nabla\mu^2\cdot \omega\nabla\omega+\mu^2\vert\nabla\omega\vert^2+G''(u)\omega^2\mu^2\right)$$
$$=\int_{\mathbb{R}^N}\left( \omega^2\vert\nabla\mu\vert^2-\mu^2\mbox{div}\left( \omega\nabla \omega\right) +\mu^2\vert\nabla\omega\vert^2+G''(u)\omega^2\mu^2\right)$$
$$=\int_{\mathbb{R}^N}\left( \omega^2\vert\nabla\mu\vert^2+\left( -\Delta\omega+G''(u)\omega \right) \omega\mu^2\right) .$$

Take $\omega:=\sqrt{G(u)}$. An easy computation shows that

$$-\Delta\omega+G''(u)\omega =-(\sqrt{G})'(u)\Delta u-(\sqrt{G})''(u)\vert\nabla u\vert^2+G''(u)\sqrt{G(u)}$$
$$=-\frac{G'(u)}{2\sqrt{G(u)}}G'(u)-(\sqrt{G})''(u)\vert\nabla u\vert^2+G''(u)\sqrt{G(u)}$$
$$=2\left( G(u)-\frac{1}{2}\vert \nabla u\vert^2 \right) (\sqrt{G})''(u).$$

Combining the above two equalities gives

$$Q(\sqrt{G(u)}\mu)=\int_{\mathbb{R}^N}\left( G(u)\vert \nabla \mu \vert^2+2\left( G(u)-\frac{1}{2}\vert \nabla u\vert^2 \right) (\sqrt{G})''(u)\sqrt{G(u)}\mu^2 \right) .$$

Finally, since $\mu$ is a $C^1$ function  with compact support in $\mathbb{R}^N$  if and only if $\sqrt{G(u)}\mu$ is also so, the proof is completed. \qed

\begin{proposition}\label{pedazodecotas}

Let $G\in C^2(\mathbb{R})$ be a nonnegative function and $u\in C^2(\mathbb{R}^N)$ be a stable nonconstant bounded entire solution of $\Delta u=G'(u)$.  Suppose that $G$ satisfies (H). Then there exist $C_1,C_2>0$ such that 
\begin{enumerate}
\item[i)] $$\int_{B_R}\left( G(u)-\frac{1}{2}\vert \nabla u\vert^2 \right)\sqrt{G(u)}\leq C_1 R^{N-2},\ \ \ \mbox{ for all }R>0.$$
\item[ii)]$$\int_{B_R}\left( G(u)-\frac{1}{2}\vert \nabla u\vert^2 \right)\leq C_2 R^{N-4/3},\ \ \ \mbox{ for all }R>0.$$
\end{enumerate}

\end{proposition}

\noindent\textbf{Proof.}

By standard regularity arguments, we can take  in Lemma \ref{caract} Lipschitz functions $\mu:\mathbb{R}^N\rightarrow \mathbb{R}$ with compact support. Fix $R>0$ and take $\mu:\mathbb{R}^N\rightarrow \mathbb{R}$ defined by
$$\mu(x):=\left\{
\begin{array}{ll}
1 & \mbox{ if } \Vert x\Vert\leq R \\ \\
\displaystyle{2-\frac{\Vert x\Vert}{R}} & \mbox{ if } R<\Vert x\Vert\leq 2R \\ \\
0 & \mbox{ if } \Vert x\Vert>2R \\ \\
\end{array}
\right.
$$

Taking into account that $0<G(u(x))\leq M$ for all $x\in\mathbb{R}^N$ (for some $M>0$), we see at once that

$$\int_{\mathbb{R}^N}G(u)\vert \nabla \mu \vert^2=\int_{B_{2R}\setminus B_R}G(u)\frac{1}{R^2}\leq \int_{B_{2R}\setminus B_R}M\frac{1}{R^2}=M(2^N-1)\vert B_1\vert R^{N-2},$$

$$ \int_{\mathbb{R}^N}\left( G(u)-\frac{1}{2}\vert \nabla u\vert^2 \right) \left( -2(\sqrt{G})''(u)\right) \sqrt{G(u)}\mu^2$$
$$\geq \int_{B_R}\left( G(u)-\frac{1}{2}\vert \nabla u\vert^2 \right) \left( -2(\sqrt{G})''(u)\right) \sqrt{G(u)}$$
$$\geq 2K \int_{B_R}\left( G(u)-\frac{1}{2}\vert \nabla u\vert^2 \right) \sqrt{G(u)}.$$

Combining these inequalities with Lemma \ref{caract} we obtain i).

To prove ii), consider the functions:
$$\alpha(x)=\left( \left( G(u(x))-\frac{1}{2}\vert \nabla u(x)\vert^2 \right)\sqrt{G(u(x)}\right)^{2/3}, \ \ \ x\in B_R,$$
$$\beta(x)=\left(\frac{G(u(x))-\displaystyle{\frac{1}{2}}\vert \nabla u(x)\vert^2}{G(u(x))}\right)^{1/3}, \ \ \ x\in B_R.$$

By i) we have
$$\Vert \alpha\Vert_{L^{3/2}(B_R)}\leq \left( C_1R^{N-2}\right)^{2/3}.$$

By Theorem \ref{modica} we see that $0\leq\beta\leq 1$ and consequently

$$\Vert \beta\Vert_{L^{3}(B_R)}\leq \Vert 1\Vert_{L^{3}(B_R)}=   \left( \vert B_1\vert R^N \right)^{1/3}.$$

Therefore, applying H{\"o}lder inequality to functions $\alpha$ and $\beta$ we conclude

$$\int_{B_R}\left( G(u)-\frac{1}{2}\vert \nabla u\vert^2 \right) =\Vert \alpha\beta\Vert_{L^1(B_R)}\leq \Vert \alpha\Vert_{L^{3/2}(B_R)}\Vert \beta\Vert_{L^{3}(B_R)}$$
$$\leq  \left( C_1R^{N-2}\right)^{2/3} \left( \vert B_1\vert R^N \right)^{1/3}= C_2 R^{N-4/3},$$

\noindent where $C_2=C_1^{2/3}\vert B_1\vert^{1/3}$, and ii) is proved. \qed

\noindent\textbf{Proof of Proposition \ref{porlomenos}.}

Applying part ii) of Proposition \ref{pedazodecotas} and Theorem \ref{modica2} we have

$$\int_{B_R}\vert \nabla u\vert^2=\int_{B_R}\left( \frac{1}{2}\vert \nabla u \vert^2 +G(u)\right)-\int_{B_R}\left(G(u)- \frac{1}{2}\vert \nabla u \vert^2 \right)$$
$$\geq  c_0R^{N-1}-C_2 R^{N-4/3},\ \ \ \mbox{ for all }R\geq 1.$$

Choosing $R_0>1$ such that $ c_0R^{N-1}-C_2 R^{N-4/3}\geq c_0R^{N-1}/2$ for every $R\geq R_0$, we complete the proof. \qed

\noindent\textbf{Proof of Proposition \ref{energiasiguales}.}

Applying Theorems \ref{modica} and \ref{modica2}, it is deduced that $$2\int_{B_R} G(u)\geq \int_{B_R}\left( G(u)+\frac{1}{2}\vert \nabla u\vert^2 \right)\geq c_0 R^{N-1},$$ for every $R\geq 1$, for some $c_0>0$. Combining this with part ii) of Proposition \ref{pedazodecotas}  we conclude

$$1\geq\frac{\displaystyle{\int_{B_R}\frac{1}{2}\vert\nabla u\vert^2}}{\displaystyle{\int_{B_R}G(u)}}=1-\, \frac{\displaystyle{\int_{B_R}\left( G(u)-\frac{1}{2}\vert \nabla u\vert^2 \right)}}{\displaystyle{\int_{B_R}G(u)}}\geq 1-\frac{C_2 R^{N-4/3}}{c_0R^{N-1}/2}=1-\frac{2C_2}{c_0R^{1/3}},$$

\noindent for every $R\geq 1$, and the proposition follows. \qed

\section{Counterexample}

\noindent\textbf{Proof of Theorem \ref{optimal}.}

Let $0<H\in C^\infty (\mathbb{R}^{N-1})$ satisfying $\int_{\mathbb{R}^{N-1}}H^2=1$. Let $k>2$ and consider an odd function $g\in C^\infty (\mathbb{R})$ satisfying $g(r)=2-r^{2-k}$ if $r\geq 1$ and $g'(r)>0$ for all $r\in \mathbb{R}$.

Define 

$$\varphi (x_1,...,x_{N}):=\frac{H(x_1,...,x_{N-1})}{\sqrt{g'(x_N)}},$$

$$\sigma(x_1,...,x_N):=g(x_N).$$

(If $N=1$, then define $\varphi(x)=1/\sqrt{g'(x)}$ and we apply the same reasoning that in the case $N>1$).

Clearly $0<\varphi \in C^\infty (\mathbb{R}^N)$ ,  $\sigma \in C^\infty (\mathbb{R}^N)$  satisfy

$$\nabla \sigma(x_1,...,x_N)=\left( 0,...,0,g'(x_N)\right),$$ 
$$(\varphi^2 \nabla \sigma)(x_1,...,x_N)=\left( 0,...,0,H^2(x_1,...,x_{N-1})\right),$$

\noindent which implies div$(\varphi^2 \nabla \sigma)=0$ in $\mathbb{R}^N.$

On the other hand, taking into account that $\sigma^2 <4$ in $\mathbb{R}^N$, and $B_R\subset \mathbb{R}^{N-1}\times (-R,R)$,  we obtain for arbitrary $R\geq 1$:

$$\int_{B_R}(\varphi \sigma)^2 dx\leq 4\int_{B_R}\varphi^2 dx\leq 4\int_{\mathbb{R}^{N-1}}H^2  \, d(x_1,...,x_{N-1})\int_{-R}^R \frac{dt}{g'(t)}$$

$$=8\int_0^R \frac{dt}{g'(t)}=8\int_0^1 \frac{dt}{g'(t)}+8\int_1^R \frac{dt}{(k-2)t^{1-k}}$$

$$=8\int_0^1 \frac{dt}{g'(t)}+\frac{8(R^k-1)}{k(k-2)}=C_1 R^k +C_2,$$

\noindent for certain $C_1>0$, $C_2\in \mathbb{R}$. Taking $C=C_1+\vert C_2\vert$, the proof is complete. \qed

\end{document}